\documentstyle[12pt]{article}
\newtheorem{lemma}{Lemma}
\newtheorem{theorem}{Theorem}

\newcommand{\QED}{{\hfill$\Box$\medskip}}
\newcommand{\k}{K\"{a}hler }

\newtheorem{corollary}{Corollary}

\begin{document}
\begin{center}
{\Large \bf Surgery and the Yamabe invariant}\\

\vspace{1 cm}

Jimmy Petean\\
Max-Planck Institut f\"ur Mathematik,\\
Bonn, Germany\\
and\\
Gabjin Yun\footnote{Supported in part by KOSEF 960-7010-201-3.}\\
Department of Mathematics\\
Myong-Ji University, Korea
\end{center}

\vspace{1 cm}

\begin{abstract}
We  study the Yamabe invariant of manifolds obtained as connected
sums along submanifolds of codimension greater than 2. In particular:
for a compact connected manifold $M$ with no metric of positive
scalar curvature, we prove that the Yamabe invariant of any manifold 
obtained by performing surgery on spheres of codimension greater 
than 2 on $M$ is not smaller than
the invariant of $M$. 
\end{abstract}

\section{Introduction}

Given a smooth compact manifold $M^n$ of dimension $n$ the {\it Yamabe
invariant} of $M$ is defined in the following way: first pick a conformal
class $\cal{C}$ of Riemannian metrics  on $M$ and let

$$Y(M,{\cal C})=\inf\limits_{g\in {\cal C}} 
\frac{\int\limits_M s_g \ \ dvol_g}{(Vol_g (M))^{\frac{n-2}{n}}}.$$

\noindent
Here $s_g$ is the scalar curvature of $g$ and  $Vol_g (M)$  is the
volume of $M$ with respect to the metric $g$. Then

$$Y(M)=\sup\limits_{\cal{C} }Y(M,\cal{C})$$

\noindent
is the Yamabe invariant of $M$
(the supremum is taken over the family of all conformal classes of metrics
on $M$). This invariant was introduced by O. Kobayashi
in \cite{Kobayashi} and  it is also 
frequently called the {\it sigma constant} 
of $M$.
Note  that $Y(M)$ is an invariant of the smooth structure of $M$.

\vspace{.5cm}

Yamabe first considered what we called $Y(M,{\cal C})$ 
in \cite{Yamabe}. He gave
a proof that $Y(M, {\cal C})$
is always achieved by a metric in ${\cal C}$. His proof contained a mistake
pointed out by Trudinger \cite {Trudinger}, and fixed in succesive steps
by Trudinger, Aubin \cite{Aubin} and Schoen \cite{Schoen2}; so proving that
in every conformal class there is a metric with constant scalar curvature.
A basic fact about the Yamabe invariant is that it is positive
if and only if the manifold admits a metric of positive scalar curvature.
The study of the invariant is then naturally divided into two cases:
when $Y(M)>0$ and when $Y(M)\leq 0$. 

\vspace{.5cm}

We will be concerned in this paper with manifolds for which the Yamabe 
invariant is non-positive. 
In this case,
every conformal class admits a unique metric of unit volume  and  
constant scalar curvature. 
That constant is precisely the Yamabe constant of the conformal
class. Hence the Yamabe invariant in this case is the supremum of the scalar
curvature over the space of metrics of constant scalar curvature
and unit volume. Also in this case, we can express the Yamabe invariant as

$$Y(M)=-\left( \inf\limits_{\cal M} {\int}_M |s_g|^{n/2} dvol_g \right)^
{2/n}, $$

\noindent
where ${\cal M}$ is the space of Riemannian metrics on $M$ (c.f. 
\cite{Anderson2}, \cite{LeBrun}). 
So the study of the 
Yamabe invariant in this case is equivalent to the study of the infimum of
this natural Riemannian functional. This is clearly not true in the positive 
case, where the infimum considered is always zero.

\vspace{.5cm}

Computations of the invariant are not easy to carry out. Nevertheless, there
has been some success in low dimensions. In dimension 4, LeBrun computed
the invariants for all compact \k surfaces with non-positive Yamabe 
invariant (see \cite{LeBrun}). The computations of the invariant for 
3-dimensional manifolds (always with non-positive invariant) are a
by-product of  Anderson's  program for the hyperbolization conjecture
(see \cite{Anderson}).
The invariants are readily computable in dimension 2 from the Gauss-Bonnet
formula. 

\vspace{.4cm}

We will prove the following:

\vspace{.5cm}

\begin{theorem} Let $M_1$, $M_2$ be compact smooth $n$-dimensional 
manifolds. Suppose that a $k$-dimensional manifold $W$ embedds into
both $M_1$ and $M_2$ with trivial normal bundle. Assume $k<n-2$.
Let $M_{1,2}^W$ be
the manifold obtained by gluing $M_1$ and $M_2$ along $W$. Then

a) If $Y(M_i)\leq 0$, then
$Y(M_{1,2}^W )\geq -[(-Y(M_1 ))^{n/2} +(-Y(M_2 ))^{n/2}]^{2/n}$.

b) If $Y(M_1 )\leq 0$ and $Y(M_2 )>0$ then $Y(M_{1,2}^W ) \geq Y(M_1)$.

\end{theorem}

\noindent
Here as usual $M_{1,2}^W$ is constructed by deletting the image of $W$ on
both $M_1$ and $M_2$  and then identifying $x\in M_1 -W$ with 
$h_2 jh_1^{-1}(x)$,
where $h_i$ is a trivialization of the normal bundle of $W$  in $M_i$ 
(i.e. a diffeomorphism between the normal bundle and a tubular
neighbourhood of $W$) and
$j$ is the inversion in ${\bf R}^{n-k} -\{ 0\}$  ($j(v)=v/{\| v\| }^2$). The 
resulting $M_{1,2}^W$ will depend on (the homotopy class of) the 
trivializations $h_i$ although we
do not make this explicit in our notation.

In particular, we have:

\begin{corollary} If $\hat{M}$ is obtained from $M$ by performing surgery
of codimension at least 3 and $Y(M)\leq 0$ then $Y(\hat{M})\geq Y(M)$.
\end{corollary}

Note that O. Kobayashi \cite{Kobayashi} has  proved the 0-dimensional
case of this result (i.e. for connected sums at points).

\vspace{.2cm}

In the next section we will carry out the main computations. We will show that
given a metric $g$ on $M$, we  can deform $g$  on the end of $M-W$ ($W$ is
an embedded submanifold of codimension at least 3) without appreciably
changing the volume or the minimum of the scalar curvature
(more precisely, for any positive $\epsilon$ we will show
that we can construct a deformation with scalar curvature bounded below
by  $s_g -\epsilon$ and volume bounded above by $Vol_g (M) +\epsilon$).
Moreover we will be able to ``choose'' how the metric looks like at the
end of $M-W$. This will be done following the constructions of Gromov-Lawson
\cite{Gromov} and Schoen-Yau \cite{Schoen}. In the last section we will
use this construction to prove the results we mentioned.

\section{The main construction}

Let $(M,g)$ be  a compact smooth Riemannian manifold  of dimension $n$.
Let $W$ be a $k$-dimensional embedded submanifold ($k\leq n-3$)
with trivial normal bundle $N$. 
Also let $S$ be the unit sphere bundle of $N$. Of course,
$M-W$ has an end diffeomorphic to $S\times (0,1)$. In this section we
will construct, for any positive (small) constant ${\epsilon}_0$, a metric 
$\hat{g}$ on $M-W$ which verifies:

a)  $\hat{g} =g$ away from the end.

b) the scalar curvature, $s_{\hat{g}}$, of ${\hat{g}}$ is greater than 
$s_g -{\epsilon}_0$ everywhere.

c) $Vol_{\hat{g}}(M-W)\leq Vol_g (M) + {\epsilon}_0$.

d)at the end $\hat{g}=h +dE^2 +dt^2$, where $h$ is any metric on
$W$ picked previously, $dE^2$ is the Euclidean metric on the sphere
$S^{n-k-1}(r)$ 
and $dt^2$ is the Euclidean metric on the line.
The positive number $r$ can be taken as small as desired.

\vspace{.5cm}

In the next section we will use this construction to prove the results 
mentioned in the introduction.

\vspace{.3cm}

We will separate our task into two steps:

\vspace{.3cm}

{\em Step 1}: 
construct $\hat{g}$ satisfying (a), (b) and (c) so that in the end it
looks like the product of the Euclidean metric on the line and the metric
$g_{\delta}$
induced by $g$ on the $\delta$-sphere bundle of $N$ (here we will be able to
pick $\delta$ as small as we want).

{\em Step 2}: 
for $\delta$  small enough find a homotopy $ds_t$, $0\leq t\leq 1$,
between
the metric $g_{\delta}$ and any metric like the one described in d) so that 
the ``total'' metric $ds_t +dt^2$ on $S\times [0,1]$
has positive scalar curvature everywhere and has volume as small as we want.
Note that we will use  this construction for metrics 
of non-positive scalar curvature,
hence this will be enough to verify condition (b).

\vspace{1cm}

Let us begin with  Step 1. The construction can be done following the
one by Gromov and Lawson \cite{Gromov} in their study of metrics with
positive scalar curvature (see also \cite{Schoen}). We will only
sketch briefly those parts where the construction in \cite{Gromov}
applies directly to our case.

First choose a positive number $\epsilon$ so that the exponential map gives
a diffeomorphism between the $\epsilon$ ball in the normal bundle and
an open neighbourhood $U_{\epsilon}$ 
of $W$ in $M$. Of course, we can pick $\epsilon$ as
small as we want. Let $N_{\epsilon}$ denote the $\epsilon$ ball in the
normal bundle, i.e. $N_{\epsilon}=\{ (x,y)\in N : \| y \| <\epsilon
\}$. Pull back the metric $g$ from $M$ to $N_{\epsilon}$ via the
exponential map. We will still call this metric on $N_{\epsilon}$ by $g$.
 
For any $0<\delta <\epsilon$ we will denote by $S^{\delta}N$ the
(n-k-1)-sphere bundle of points in $N$ of norm $\delta$ and we will
call $g_{\delta}$ the restriction of $g$ to  $S^{\delta}N$. Also, for
any point $(x,y)\in N$, let $r(x,y)=\| y\|$ denote the  distance (measured
by $g$) from $(x,y)$ to $W$. 
 
\vspace{.2cm}

Now pick any $r_1$  so that $0<r_1 <\epsilon$. In the plane with coordinates
$(t,r)$ consider a smooth curve $\gamma$
which moves vertically from $(0,\epsilon )$
to $(0,r_1 )$, then is the graph of a  strictly decreasing 
function $h$ joinning
the points $(0,r_1 )$ and $(t_f ,r_f )$, and then stays on the horizontal
line $r=r_f$ (here $t_f$ is a positive number and $0<r_f <r_1$).   

\vspace{.2cm}

Let $M^{\gamma} =\{ ((x,y), t)\in N_{\epsilon} \times {\bf R} : (r(x,y),t)
\in \gamma \}$. Let $g^{\gamma}$ denote the restriction to $M^{\gamma}$
of the product metric $g+dt^2$ ($dt^2$ denotes the Euclidean metric on the
line). To complete Step  1 we need to show that we can pick $\gamma$ so
that the metric $g^{\gamma}$ has the required restrictions on scalar
curvature and volume (from the shape of $\gamma$ one can immediately see
that $g^{\gamma}$ will verify condition (a) and will have the required
form at the end). 

Let us study the volume of $g^{\gamma}$. More precisely, we will describe
what  conditions must verify $\gamma$ in order  that the volume
of $M^{\gamma}$ is bounded by some chosen  small positive constant $u$.

\vspace{.2cm}

The horizontal piece of $\gamma$ can be taken arbitrarily small. So we
only need to consider the volume of the part of $\gamma$ corresponding
to the graph of $h$.

We will then restrict our attention to the
piece of $\gamma$ given by the graph of $h$. But to avoid introducing
new notation we will keep calling  $\gamma$  the smaller curve
and $g^{\gamma}$ the restriction of the metric to this piece.
To estimate the volume of $g^{\gamma}$ we will  compare its volume element
to the one of the original metric $g$.
Of  course, $dvol_{g^{\gamma}} = dvol_{g^{\gamma}}(e_1 ,...,e_n )dvol_g$;
where $e_1 ,...,e_n$ form a $g$-orthonormal basis (at certain point).

Let $h^{-1}:(r_f ,r_1 )\rightarrow (0,t_f )$ be the inverse function of $h$.
Consider the map

$$F:N_{r_1} -\overline{N_{r_f}}\rightarrow M^{\gamma}$$

\noindent
given by $F(x,y)=(x,y,h^{-1}(r(x,y)))$. $F$ is clearly  a diffeomorphism
(onto the part of  $\gamma$ corresponding to the graph of $h$)
and we will use it to bring the metric $g^{\gamma}$ to 
$N_{r_1}-\overline{N_{r_f}}$.

The function $r:N_{r_1}-\overline{N_{r_f}} \rightarrow {\bf R}$  
is a smooth  submersion. So,  at  
any point $q=(x,y)\in N_{r_1}-\overline{N_{r_f}}$, the kernel  of $r_{* q}$
has dimension $n-1$. 
Let $e_2 ,...,e_n$ be a $g$-orthonormal basis of Kernel($r_{* q}$). And
let $e_1$ be a  unitary vector orthogonal to that space.

\vspace{.3cm}

We have to compute $(dvol_{g^{\gamma}})_{F(q)} 
(F_{* q}e_1 ,...,F_{* q}e_n)$.
Now for $i=2,..,n$ we have $F_{* q}e_i =(e_i ,0)$ as elements of
$T_q M \bigoplus T_{F(q)}{\bf R}$. Moreover,

$$(dvol_{g^{\gamma}})_{F(q)} (F_{* q}e_1 ,...,F_{* q}e_n )=
\sqrt{\det (g^{\gamma}(F_* e_i ,F_* e_j )_{ij})}.$$

It follows from the previous considerations that the (n-1)$\times$(n-1) matrix 
obtained by deletting the first row and column of 
$g^{\gamma}(F_* e_i ,F_* e_j )_{ij}$
is the identity.
But also,

$$g^{\gamma}(F_{*}e_i ,F_{*}e_1 )={\delta}_{i1}+
dt^2 ((h^{-1} r)_{*} e_i , (h^{-1} r)_{*}e_1 ).$$

\noindent
Hence the matrix is diagonal, having 1 in all but the first diagonal
entry and $1+dt^2 ((h^{-1} r)_{*} e_i , (h^{-1} r)_{*}e_1 )$ in
the first.

We need to compute, or estimate,   $(h^{-1} r)_{*}e_1 $.   

\vspace{.5cm}

{\it Claim}: Pick any small positive constant $u$. 
If $r_1$ is small enough then at any $q$ in $N_{r_1}-W$
we have $|r_{* q}e_1 -1|\leq u$. Here we are identifying the
tangent space of the  real line with the real numbers, and we
are picking the $e_1$ `` with the right sign''.

\vspace{.3cm}

To prove
the claim, let $(x_1 ,...,x_k )$ be coordinates for $W$ on an open
set $U$. Pick an orthonormal frame for $N$ over $U$ and let
$(y_1 ,...,y_{n-k})$ be the induced coordinates on the fibers.
Hence $(r(x,y))^2 =y_1^2 +...+y_{n-k}^2$.

In $U\times D^{n-k}  (r_1 )$ it  is easy to find a $g$-orthonormal
basis for the kernel of
$r_{*}$. At a point $(x,y)$ let $\frac{\partial}{\partial r}$  be the
unitary vector tangent to the geodesic line from  0 to $y$. Let
$v_2 ,...v_{n-k}$ be an orthonormal basis of the tangent space to the
sphere of radius $r(x,y)$. Let $w_{n-k+1},...,w_n $ be a $g$-orthonormal
basis of the tangent space to $U$ at $x$.   
Note that at $q\in N_{r_1}$, 
$$g(v_i ,w_j )\in o(r_1 ) , \ \ \ \ 
 g(\partial /\partial r, w_j )\in o(r_1 )$$

\vspace{.2cm}
\noindent
and since we are in
geodesic coordinates $g(\frac{\partial}{\partial r}, v_i )=0$.
Applying the Gram-Schmidt process to the basis $v_2 ,...,v_{n-k},w_{n-k+1},
...,w_{n-1}, \frac{\partial}{\partial r}$ we get an orthonormal basis.
The first $n-1$ vectors will be an o.n. basis for the kernel of $r_{*}$.
The last vector, which is the vector $e_1$ in the claim, 
will differ from $ \frac{\partial}{\partial r}$ by
$o(r_1 )$. But it is clear that $r_{*}( \frac{\partial}{\partial r})=
 \frac{\partial}{\partial t}$ and the claim follows.

\vspace{.4cm}

Hence, with the notation we were using before the claim,

$$g^{\gamma} (F_{* q}e_1 , F_{* q}e_1)=1+(1+o(u))^2 
\left( \frac{\partial h^{-1}}{\partial t} (r(q)) \right)^2$$.

\vspace{.2cm}

Finally, 

$$dvol_{g^{\gamma}}(e_1 ,...,e_n )=\sqrt{1+(1+o(u))^2 
\left( \frac{\partial h^{-1}}{\partial t}
 (r(q))\right) ^2}\leq $$
$$\leq 1-(1+o(u)) \frac{\partial h^{-1}}{\partial t}(r(q))$$

\noindent
and we have


\begin{eqnarray}
Vol_{g^{\gamma}}(N_{r_1}-N_{r_f})\leq Vol_g (N_{r_1} )-
(1+o(u)) \int\limits_{N_{r_1} -N_{r_f}} 
\frac{\partial h^{-1}}{\partial t} (r(q)) \  dvol_g .
\end{eqnarray}

To estimate the last integral we will need the following result, whose 
proof is elementary.

\begin{lemma} Given $\delta >0$ small enough, there exists a 
constant $\cal{K}$, depending on $g$ and $\delta$,
such that for all $0<r_a <r_b <\delta$,

$$Vol_g (N_{r_b} -N_{r_a})\leq (r_b -r_a )\cal{K}.$$

\end{lemma}

Now we can see what are the conditions we need for $\gamma$.
Suppose we want to make the volume of the portion of $M-W$  with
the newly defined metric less than ${\epsilon}_0$. Picking $r_1$ small
enough we can make the first term in $(1)$ less than ${\epsilon}_0 /2$.
It follows from the previous lemma  that the second term in $(1)$
is bounded by $2(1+o(u)){\cal K}\  t_f$. 
As we let $r_1$ tend to 0, the
first 3 factors in the previous product are bounded independently of
$r_1$ (with a little effort we  could prove that the product actually
tends to 0, but we will not need this). 

\vspace{.3cm}

So, to keep the volume under control, 
we just need to construct $\gamma$ so that as $r_1$ tends
to 0, $t_f$ also tends to 0.

\vspace{.2cm}

We will now show that we can construct such a curve $\gamma$ which satisties
also our requirements about the scalar curvature. More precisely, given
any small positive $r_1$ we will construct $\gamma$ so that $t_f \leq 7r_1$
and $s_{g^{\gamma}}\geq s_g -{\epsilon}_0$ (where ${\epsilon}_0$  is a
positive number we fixed previously).

\vspace{.8cm}

Let $\theta$ be the angle between the $r$-axis
and the tangent to $\gamma$ (at some point).  
Let $\kappa$ be the principal curvature  of $\gamma$.
We have the following formula for the scalar curvature of $g^{\gamma}$
(it is computed in \cite{Gromov}),

$$s_{g^{\gamma}}=s_g +O(1){\sin }^2\theta +\frac{(n-k-1)(n-k-2)}{r^2} 
{\sin }^2\theta -(n-k-1)\frac{\kappa}{r} \sin \theta.$$

\noindent
The function $O(1)$ is bounded independently of $r$, $\theta$ and $\kappa$. 
Let
${\cal A}$ be an upper bound for it.

To construct $\gamma$ we will be describing how we make $\theta$ vary from
0 to $\pi /2$. Let $s$ be the arclength parameter. Let $\kappa (s)$ 
be the curvature
of $\gamma$ at (length) $s$. The change of $\theta$  between $s_1$ and
$s_2$ is given by the integral of $\kappa (s)$ 
(between $s_1$ and $s_2$, of course).

\vspace{.3cm}

We will describe  the function $\kappa (s)$.
First pick ${\theta}_0$ so that


\begin{eqnarray}
\left({\cal A}
+\frac{4n}{r_1^2 }\right) \sin {\theta}_0 <{\epsilon}_0. 
\end{eqnarray}

Bend $\gamma$ changing $\theta$ from 0 to ${\theta}_0$ by a curve of length
$r_1$/2. We can do this keeping $\kappa (s)$ bounded by $2/r_1$, 
and  so it is easy to
see from (2) that we will have 
$s_{g^{\gamma}}\geq s_g -{\epsilon}_0$. We will continue
on the straight line of angle ${\theta}_0$ till we get to a point
with $r=r_2 <(\sin {\theta}_0 )/{\cal A}$.
We will now make a new
bend of length $r_2$/2. In order to keep the desired inequality for the 
scalar curvature it is enough to keep $\kappa (s)$ 
bounded by $(\sin {\theta}_0 )$/2$r_2$;
hence we will increase $\theta$  by about $(\sin {\theta}_0 )$/4. Then we will
continue again in a straight line (of angle ${\theta}_1$) for an 
arbitrarily small period  of time.
We will repeat this process until we got  to a horizontal line 
($\theta = \pi$/2). 
Note that the amount of $\theta$ we can increase at every bend
is about $(\sin {\theta}_{n-1})$/4 (where ${\theta}_{n-1}$ is the angle
previous to the new bending). It is then clear that we need to repeat the
process only a finite number of times. Now pick $i$ so that $\sin {\theta}_i$
is about 1/2. Then we will clearly finish the process at the step $i+10$.
Hence the  length of $\gamma$ in this last part of the process is bounded
by $5 r_1$.
On the other hand all the previous bends will verify that $r_{k-1} -r_k
\geq (1/4) r_{k-1}$. 
This implies that the length of the first bends 
(and straight lines) are bounded by twice the
decrease in $r$. Hence the total length of ${\gamma}$ up to the
$i$-th step is bounded by 2$r_1$. We can then construct the whole $\gamma$
of length bounded by 7$r_1$. Therefore $t_f \leq 7r_1$. This completes 
Step 1 of the construction.

\vspace{.6cm}

{\it Step} 2: We have already deform $g$ so that in the end we have the
product metric $g_{\delta} +dt^2$. We now want to find a homotopy 
between $g_{\delta}$ and some canonical metrics. We will have a family 
$g^t$, $a\leq t\leq b$ of metrics of positive scalar curvature on a compact
manifold $X$ and
we need to know conditions
under which the metric $g^t +dt^2$, on $X\times [a,b]$,
has positive  scalar curvature and small volume.

\vspace{.3cm}

The following two lemmas are the tools we need to study this problem. 
The first one will give us a bound on the total
volume of  a homotopy and is completely elementary,

\begin{lemma} Let $g^t$, $a\leq t\leq b$, be a family of  metrics on the
compact manifold $X$. Let $dt^2$ be  the Euclidean metric on the interval
$(a,b)$ and consider the metric $G=g^t +dt^2$ on $X\times (a,b)$. We have,
$$Vol_G (X\times (a,b))\leq \sup\limits_{a\leq t\leq b} \{ Vol_{g^t} (X)\}
(b-a).$$
\end{lemma}

We will be concerned with the case in which all the metrics $g^t$ have strictly
positive scalar  curvature. One can then see (\cite{Gromov}, \cite{Schoen}) 
that
doing the homotopy slowly in time (stretching the interval) the  scalar
curvature of the total metric is also positive. We will now write down
explicitly how much we need to stretch the interval  in order to get
positive scalar curvature. 

First we need to introduce some notation. Suppose  we have a homotopy,
$H=\{ g^t \} $, $0\leq t\leq 1$, of positive  scalar
curvature metrics on a compact manifold $X$. Let 

$$s(H)=\min\limits_{(x,t)} \{ s_{g^t} (x)\} .$$

Now consider a  system of coordinates $\bar{x} =(x_1 ,...,x_n )$ 
in some open neighbourhood in $X$.
Let $G_a$ be the metric $ g^t +a^2 dt$. 
Let $g_{ik} (x,t)=g^t (\partial /\partial x_i ,\partial /\partial x_k )$
and let $g^{ik}$ denote as usual the coefficients of the inverse of the
matrix $(g_{ik} )$. 
By a straightforward computation
we get the following formula,

$$s_{G_a} (x,t)=  s_{g^t}(x)+
\frac{1}{4a^2}\left( g^{ik} 
\frac{\partial g_{jk}}{\partial t} g^{jp} \frac{\partial g_{ip}}{\partial t}
- g^{ik} \frac{\partial g_{ik}}{\partial t}  g^{jp}
\frac{\partial g_{jp}}{\partial t} \right)(x,t) $$
$$\ \ \ \ \ \ \ \ \ \ \ \ \ \ \ \ \ 
-\frac{1}{2a^2} \left(
g^{ik}\frac{{\partial}^2 g_{ik}}{\partial t^2}
+\frac{\partial}{\partial t} (g^{ik}
\frac{\partial g_{ik}}{\partial t}) \right)(x,t),
$$
where $s_{G_a}$ and $s_{g^t}$ denote the scalar curvatures of the
metrics $G_a$and $g^t$, respectively.

We will call ${\cal B}(x,t)$ the expression multiplied by 1/$a^2$ in the
previous formula and ${\cal B}(H)$ the maximum of ${\cal B}(x,t)$ on 
$X\times [0,1]$. Note that ${\cal B}(x,t)$ does not depend on the system
of coordinates.
We have (compare to \cite[Lemma 3]{Gromov}

\begin{lemma} Let $H=g^t $, $0\leq t\leq 1$, be a family of metrics of 
strictly
positive scalar curvature on the compact manifold $X$. The metric
$\hat{G}_a =g^{t/a} +dt^2$ 
on $X\times [0,a]$ has positive  scalar curvature  if
$a>\sqrt{{\cal B}(H)/s(H)}$.
\end{lemma}

{\it Proof}: Consider  the map between $X\times [0,1]$ and
$X\times [0,a]$ given by multiplication by $a$ in the $t$-coordinate.
Pulling back the metric $\hat{G}_a $ via this  diffeomorphism we get
the metric $G_a =g^t +a^2 dt^2$ on $X\times [0,1]$.
The lemma now follows directly from
our formula  for the scalar curvature of $G_a$.

\QED

Now let us come back to our problem. We want to study the metric $g$
close to $W$. Find a coordinate system of the form 
$(x_1, \cdots, x_k, y_{k+1}, \cdots, y_n)$. Where the $x_i$ are coordinates for $W$ and the $y_i$ are coordinates
on the normal bundle induced by an orthonormal  frame  (make the identification
via the exponential map). As before  $r(x,y)$ will denote the distance between
the point $(x,y)$ and $W$. So $r^2 =\Sigma y_i^2$.

\vspace{.3cm}

The metric $g$ can be expressed in these coordinates in the form,

$$g=\Sigma (g_{ij} (x,0)+o(r))\  dx_i dx_j +\Sigma o(r)  \ 
dx_i dy_{\alpha}+$$

$$\ \ \ \ \ \ \ \ \ \ \ \ \ \ \ \ +\Sigma ({\delta}_{\alpha \beta} 
+o(r^2 ))\ dy_{\alpha} dy_{\beta} .$$

We will first homotope $g_{\delta}$ to the product metric 
$g_{\delta}^N =g|_W +dE^2  (\delta )$  
on $W\times S^{n-k-1}$ ($dE^2 (\delta )$ is the standard
Euclidean metric on the $\delta$-sphere). It is easy to see from our
local formula that the metrics $g_{\delta}$ and $g_{\delta}^N$ will
differ only by low order terms. Consider the simplest homotopy
$g_{\delta}^t =tg_{\delta}^N +(1-t)g_{\delta}$. 
In local coordinates $g_{\delta}^t$ will look like $g_{\delta}^N$ plus some
terms of order $o(r)$ multiplied by $t$.
If we call $H^1$ this family
of metrics we have then that ${\cal B}(H^1 )$ is bounded independently
of $\delta$ (for $\delta$ smaller than some fixed small ${\delta}_0$).
It is also clear that $s(H^1 )$ is bounded below by some positive multiple
of $1/{\delta}^2$. We also have that the volume of the  metrics 
$g_{\delta}^t$ go
to 0 as we take $\delta$ small. It follows from the previous lemmas
that given any ${\epsilon}_0 >0$ for any  $\delta$ small enough
the metric $g_{\delta}^t +dt^2$ on $W\times S^{n-k-1}\times [0,1]$ will have
positive scalar curvature everywhere and volume less than ${\epsilon}_0$.

\vspace{.2cm}

Hence we can assume now that at the end we have the product metric, 
$g_{\delta}^N$,
of the Euclidean metrics (on the $\delta$-sphere and on the line) and
the restriction of $g$ to $W$.

\vspace{.2cm}

Finally, we need to show that we can change the metric $g|_W$ to any
other metric $h$ on $W$. Consider again the homotopy $g^t =
t(h+dE^2 )+(1-t)g_{\delta}^N$. Of course, the homotopy only works on the 
$W$-part, while we always have the Euclidean metric on the $\delta$-sphere.
Let us  call $H^2$ this homotopy. It is again easy to see that
${\cal B}(H^2 )$ is bounded independently of $\delta$ (it is actually
independent of $\delta$), while $s(H^2 )$
is bounded below by some positive multiple of $1/{\delta}^2$.  
The  volumes of the  metrics $g^t$ tend to 0. These observations
and the previous lemmas show that if $\delta$ is small enough 
the metric $g^t +dt^2$ on  $W\times S^{n-k-1}\times [0,1]$ has positive
scalar curvature and  volume less than ${\epsilon}_0$.

\vspace{.3cm} 

We have therefore finished Step 2 and so our construction.

\section{Proof of Theorem 1}

In this section we will use the previous construction to prove the results
mentioned in the introduction.  We will use the following well known result
(see \cite{Kobayashi}),

\begin{lemma} Let $\cal{C}$ be a conformal class of Riemannian metrics
on the compact $n$-dimensional manifold $M$. Assume $n\geq 3$ and 
$Y(M,{\cal C} )\leq 0$. Let $g$ be any metric in ${\cal C}$. Then,

$$Y(M,{\cal C} )\geq min(s_g ) (Vol_g (M))^{2/n}.$$

\end{lemma}

\vspace{.3cm} 

{\it Proof of Theorem 1}:
We will use the same notations as in the introduction.
To begin, fix any positive $\epsilon$. Pick conformal classes ${\cal{C}}_i$ in
$M_i$ so that $Y(M_i ,{\cal{C}}_i )\geq Y(M_i )-\epsilon$ (for $i=1,2$).
Let $g_i$ be the metric of volume ${\lambda}_i$ in ${\cal{C}}_i$ 
realizing $Y(M_i,{{\cal C}}_i )$. 
Renormalize the metrics so that ${\lambda}_1 +{\lambda}_2 =1$.

Now construct metrics $\hat{g_i}$ in $M_i -W$, as in the previous
section (using the $\epsilon$ we have just picked),  
so that they coincide in the
ends. Then we can glue them together 
to obtain a metric $\hat{g}$ on $M_{1,2}^W$.
The volume of $\hat{g}$ satisfies

$$ Vol_{\hat{g}} (M_{1,2}^W )\leq 1+2\epsilon$$

Note also that the scalar curvature of $g_i$ is constant equal to
$$s_{g_i}=Y(M_i,{\cal{C}}_i )({\lambda}_i )^{-2/n}.$$

Hence,

$$s_{\hat{g}}\geq min \left\{ 
\frac{Y(M_1 ,{{\cal C}}_1 )}{({\lambda}_1 )^{2/n}}
-\epsilon ,
\frac{Y(M_2 ,{{\cal C}}_2 )}{({\lambda}_2 )^{2/n}} -\epsilon \right\} $$

\noindent
everywhere. So, from the previous lemma, we get:

$$Y(M_{1,2}^W , {{\cal C}}_{\hat{g}} )\geq (1+2\epsilon )^{2/n} 
min \left\{ \frac{Y(M_1 ,{{\cal C}}_1 )}{({\lambda}_1 )^{2/n}}
-\epsilon ,
\frac{Y(M_2 ,{{\cal C}}_2 )}{({\lambda}^2 )_{2/n}} -\epsilon \right\}$$

To simplify the notation, we will write $a_i =Y(M_i, {{\cal C}}_i )$. 
Assume first that
both $a_1$ and $a_2$ are strictly negative. We can rescale the metrics  
so that ${\lambda}_1$ is any number in (0,1), and 
${\lambda}_2=1-{\lambda}_1$. The minimum appearing in the previous 
formula will be the greatest when the two numbers are equal. This
happens when

$${\lambda}_i =\frac{|a_i |^{n/2}}{|a_1 |^{n/2}+|a_2 |^{n/2}}$$

\noindent
Then we get 

$$min \left\{ \frac{Y(M_1 ,{{\cal C}}_1 )}{({\lambda}_1 )^{2/n}}
-\epsilon ,
\frac{Y(M_2 ,{{\cal C}}_2 )}{({\lambda}_2 )^{2/n}} -\epsilon \right\}=
-(|a_1 |^{n/2}+|a_2 |^{n/2})^{2/n}-\epsilon$$

\noindent
and so,

$$Y(M_{1,2}^W ,{\cal C}_{\hat{g}})
\geq -(1+2\epsilon)^{2/n} \left[ \left ((-Y(M_1 )+\epsilon )^{n/2}+
(-Y(M_2 )+\epsilon )^{n/2}\right) ^{2/n}+\epsilon \right].$$

\noindent
Since $
Y(M_{1,2}^W )=sup_{{\cal C}} \{ (Y(M_{1,2}^W , \cal{C})\}$ and 
$\epsilon$ was arbitrary, we get

$$Y(M_{1,2}^W )\geq - [(-Y(M_1 ))^{n/2}+
(-Y(M_2 ))^{n/2}]^{2/n},$$

\noindent
finishing the proof of Theorem 1 in this case.
If $a_i =0$ the result follows easily making ${\lambda}_i $ tend to 0.
Part b) of Theorem 1 follows in the same way
since the fact that $Y(M_2 )>0$ implies that $M_2$ 
admits a scalar flat metric. We  have therefore
finished the proof of Theorem 1.

\QED

\noindent
{\bf Acknowledgements:} The first author would like to thank the
hospitality of the directors and staff of the Max-Planck Institut,
where he stayed during the preparation of this work. He is also very
grateful to Claude LeBrun for very helpful observations on previous
versions of this paper and to Vyacheslav Krushkal for many explanations
about surgery techniques. The second author would like to express his gratitude to M. T. Anderson for advice.

\vspace{1cm}

\noindent
{\bf Authors' addresses:}

\vspace{.5cm}

\noindent
Jimmy Petean\\
Max-Planck-Institut f\"{u}r Mathematik\\
Bonn, Germany.\\
email: petean@mpim-bonn.mpg.de

\vspace{.5cm}

\noindent
Gabjin Yun\\
Department of Mathematics\\
Myong Ji University\\
Yong In, Kyung Ki, Korea 449-728\\
email: gabjin@wh.myongji.ac.kr

\end{document}